\documentclass[11pt]{article}
\usepackage{times}
\usepackage{amsmath,amsfonts,amssymb,latexsym,epsfig}
\usepackage{color,epsf,float}
\usepackage{graphicx,stix}
\usepackage[toc]{appendix}

\usepackage{alltt}
\usepackage{url}
\usepackage{a4}

\newcommand{\e}{\ensuremath{\mathrm{e}}}
\newcommand{\ad}{\ensuremath{\mathrm{ad}}}

\begin{document}

\title{A unifying framework for perturbative exponential factorizations}

\author{Ana Arnal\thanks{Email: \texttt{ana.arnal@uji.es}}
	\and
	Fernando Casas\thanks{Email: \texttt{Fernando.Casas@uji.es } }\and Cristina Chiralt\thanks{Email: \texttt{chiralt@uji.es }} \and José Angel Oteo\thanks{Email: \texttt{oteo@uv.es}} }

\date{}

\maketitle

\noindent {\footnotesize $^{\ast}$, $\dag$,$\ddag$ Departament de Matem\`atiques and IMAC, Universitat Jaume I, 12071 Castell\'on, Spain, \\
$\S$ Departament de F\'{\i}sica Te\`{o}rica, Universitat de Val\`{e}ncia and Institute for Integrative Systems Biology (I2SysBio),  46100-Burjassot, Valencia, Spain }

\begin{abstract}
	
We propose a framework where Fer and Wilcox expansions for the solution of differential equations
are derived from two particular choices for the initial transformation that seeds the product expansion. In this scheme intermediate expansions can also be envisaged. Recurrence formulas are developed. 
A new lower bound for the convergence of the Wilcox expansion is provided as well as some applications of the results.  In particular, two examples are worked out up to high order of approximation to illustrate the behavior of the Wilcox expansion.
	\vspace*{1cm}

\noindent {\bf Keywords:} Sequences of linear transformations; Wilcox expansion; Fer expansion; Zassenhaus formula; Bellman problem.

\end{abstract}

\section{Introduction}

	Linear differential equations of the form
\begin{equation}   \label{ccv.1}
\dot{x} \equiv \frac{dx}{dt} = A(t) x,  \qquad x(0) = x_0 \in \mathbb{C}^d,
\end{equation}
with $A(t)$ a $d \times d$ matrix whose entries are integrable functions of $t$, are ever-present in many branches of science, 
{the fundamental evolution equation of Quantum Mechanics, the Schr\"odinger equation, being a particular case. In consequence, solving
	Eq. (\ref{ccv.1}) is of the greatest importance.}
In spite of their apparent simplicity, however, they are seldom solvable in terms of elementary functions, and so different procedures have been
proposed along the years to render approximate solutions. These are specially useful in the analytical treatment of perturbative problems, such
as those arising in the time evolution of quantum systems \cite{wilcox67eoa}, control theory or problems where time-ordered products
are involved \cite{oteo00fto}. Among them, exponential perturbative expansions have received a great deal of attention, due
to some remarkable properties they possess. In particular, if Eq. (\ref{ccv.1}) is defined in a Lie group, the
approximations they furnish also evolve in the same Lie group. As a consequence, important qualitative properties of the exact solution are also 
preserved by the approximations. Thus, if Eq. (\ref{ccv.1}) represents the time-dependent Schr\"odinger equation, then the
approximate evolution operator is still unitary, and as a consequence the total sum of (approximate) transition probabilities is the unity, no matter where the
expansion is truncated. {There are in fact many physical problems (in non-linear mechanics, optical spectroscopy, magnetic resonance, etc.)
	involving periodic fast-oscillating external fields that are also modeled by Eq. (\ref{ccv.1}), with $A(t)$ periodic. In that case, especially
	tailored expansions incorporating the well-known Floquet theorem \cite{hale80ode}, such as the average Hamiltonian theory 
	\cite{maricq82aoa} and the Floquet--Magnus
	expansion \cite{casas01fte,mananga16otf} have also been proposed.}

{When dealing with the general problem (\ref{ccv.1})}, 
one of the most widely used exponential  approximations corresponds to the Magnus expansion \cite{magnus54ote},
\begin{equation}
x(t) = \e^{\Omega(t)} x_0,
\end{equation}
where $\Omega$ is an infinite series,
\begin{equation}
\Omega(t) = \sum_{k=1}^{\infty} \Omega_k(t)
\end{equation}
whose terms are linear combinations of time-ordered integrals of nested commutators of $A$ evaluated at different times (see
\cite{blanes09tme} for a review, including applications to several physical and mathematical problems). What is more interesting for our purposes here is that
this expansion can be related with a coordinate transformation $x \longmapsto X$ rendering the original system (\ref{ccv.1}) into the trivial equation
\begin{equation} \label{trivial}
\frac{dX}{dt} = 0,
\end{equation}
with {the \textit{static}} solution $X(t) = X(0) = x_0$, and that the transformation is given precisely by $x(t) = \exp(\Omega(t)) X(t)$      
\cite{casas19cco}. 

In contrast to the Magnus expansion,  {  the Floquet-Magnus expansion gets the solution  with two exponential transformations when $A(t)$
	is periodic, whereas  }
other exponential perturbative expansions are based on infinite {product }factorizations of $x(t)$,
\begin{equation}  \label{facto1}
x(t) = \e^{\Omega_1(t)} \, \e^{\Omega_2(t)} \cdots \e^{\Omega_n(t)} \cdots \, x(0),
\end{equation}
such as those proposed by Fer and Wilcox. 
In fact, as pointed out in \cite{blanes98maf}, both expansions have a    
curious history worth to be sketched. It was Fer who proposed the expansion that bears his name in \cite{fer58rdl}, although
he never applied it to solve any specific problem. Bellman reviewed this paper in the \textit{Mathematical Reviews}\footnote{MR0104009.}, 
and even proposed the expansion as an exercise
in \cite{bellman70itm}. Nevertheless, Wilcox identified it in \cite{wilcox67eoa} with an alternative factorization, Eq. (\ref{facto1}), which was indeed a different
and new type of expansion. From them on, their historical trajectories move apart. Thus, Fer expansion was rediscovered by Iserles 
\cite{iserles84slo} as a tool for the numerical integration of linear differential equations and later on used in quantum mechanics \cite{klarsfeld89eip} 
and solid-state
nuclear magnetic resonance \cite{mananga16otf2}, but also as a Lie-group integrator \cite{casas96ffa,zanna99car,iserles00lgm}. 
On the other hand, Wilcox expansion has been rediscovered several
times in the literature, in particular in \cite{huillet87lac} in the context of nonlinear control systems, and in \cite{zagury10ueo} as a general tool for approximating the
time evolution operator in quantum mechanics.

The first goal in this work is to recast both infinite product expansions within a unifying framework. This is done by considering, instead of
	just \emph{one} exponential transformation as in the case of the Magnus expansion, a \emph{sequence} of such transformations,
	$\e^{\Omega_1(t)}, \ldots, \e^{\Omega_k(t)}, \ldots$, chosen to satisfy certain requirements. To be more specific, suppose one 
	replaces $A(t)$ in Eq. (\ref{ccv.1}) by $\lambda A(t)$, where $\lambda > 0$ is a parameter. Then, if the transformations are chosen so that each
	$\Omega_k(t)$ is proportional to $\lambda^k$, we recover the Wilcox expansion, whereas we end up with the Fer expansion when each 
	$\Omega_k(t)$ is an infinite series in $\lambda$ whose first term is proportional to $\lambda^{2^{k-1}}$.

We also show that further alternative descriptions yield new factorizations. This additional degree of freedom can be indeed used 
	to deal better with the features of the matrix $A(t)$ as  in Floquet-Magnus, when $A(t)$ is periodic.

One might then consider this sequence of linear transformations as a generalization of the concept of picture in Quantum Mechanics
	when Eq. (\ref{ccv.1}) refers to the Schr\"odinger equation.


Our second goal consists in obtaining, on the basis of this framework,  new results concerning Wilcox expansion. Thus, we develop a recursive procedure to obtain every order of approximation in terms of nested commutators, as well as a convergence radius bound. 
	We also establish a formal connection of the Wilcox expansion with the   Zassenhaus formula \cite{magnus54ote,casas12eco}. 

Eventually, the important problem of expanding the exponential $\exp(A + 
\varepsilon B)$ for $\varepsilon > 0$ and $A$, $B$, two generic {non-commuting} operators will be addressed, and two applications of the results given.

	\section{A sequence of transformations{: the general case}}
	
	Given the initial value problem (\ref{ccv.1}), let us consider a linear change of variables $x \longmapsto X_1$ of the form
	\begin{equation}   \label{gt.1}
	x(t) = \e^{\Omega_1(t)} X_1(t),  \qquad \Omega_1(0) = 0
	\end{equation}
	transforming the original system into
	\begin{equation}   \label{gt.2}
	\frac{d X_1}{dt} = B_1(t) X_1.
	\end{equation}
	For the time being, the generator $\Omega_1(t)$ of the transformation is not specified.
	Then, $B_1(t)$ can be expressed in terms of $A(t)$ and $\Omega_1(t)$ as follows. First, by inserting Eq. (\ref{gt.1}) into Eq. 
	(\ref{ccv.1}), and
	taking Eq. (\ref{gt.2}) into account, one can obtain
	\begin{equation}
	\frac{d }{dt} \exp(\Omega_1) = A(t) \exp(\Omega_1) - \exp(\Omega_1) B_1(t),
	\end{equation}   
	whence
	\begin{equation}
	B_1(t) = \e^{-\Omega_1} A(t) \e^{\Omega_1} - \e^{-\Omega_1} \frac{d }{dt} \e^{\Omega_1}. 
	\end{equation}
	The derivative of the matrix exponential can be written as \cite{blanes09tme}
	\begin{equation}
	\frac{d }{dt} \exp(\Omega_1(t)) = d \exp_{\Omega_1(t)} (\dot{\Omega}_1(t)) \, \exp(\Omega_1),
	\end{equation}
	where the symbol $d \exp_{\Omega}(C)$ stands for the (everywhere convergent) power series
	\begin{equation}
	d \exp_{\Omega}(C) = \sum_{k=0}^{\infty} \frac{1}{(k+1)!} \ad_{\Omega}^k (C) \equiv \frac{ \exp(\ad_{\Omega}) - I}{\ad_{\Omega}} (C).
	\end{equation}
	Here $\ad_{\Omega}^0 C = C$, $\ad_{\Omega}^k C = [\Omega, \ad_{\Omega}^{k-1} C]$,  and $[\Omega,C]$ denotes the usual commutator. 
	Therefore
	\begin{equation}  \label{b1}
	B_1(t) = \e^{-\Omega_1} \left( A(t) -   d \exp_{\Omega_1(t)} (\dot{\Omega}_1(t)) \right) \e^{\Omega_1}  \equiv  \e^{-\ad_{\Omega_1}} (B_0 - G_1)
	\end{equation}   
	where 
	\begin{equation}  \label{b2}
	B_0(t) \equiv A(t), \qquad G_1(t) \equiv d \exp_{\Omega_1(t)} (\dot{\Omega}_1(t)) 
	\end{equation}   
	and   
	\begin{equation}
	\e^{-\ad_{\Omega_1}} F = \sum_{k \ge 0} \frac{(-1)^k}{k!} \ad_{\Omega_1}^k F = \e^{-\Omega_1} F \, \e^{\Omega_1}.
	\end{equation}  
	Of course, nothing prevents us from  repeating the whole procedure above and introduce a second transformation to Eq. (\ref{gt.2}) of the form
	\begin{equation}
	X_1(t) = \e^{\Omega_2(t)} X_2(t),  \qquad \Omega_2(0) = 0,
	\end{equation}
	so that the new variables verify
	\begin{equation}
	\frac{d X_2}{dt} = B_2(t) X_2.
	\end{equation}
	In general, for the $n$-th such linear transformation
	\begin{equation}
	X_{n-1}(t) = \e^{\Omega_n(t)} X_n(t),  \qquad \Omega_n(0) = 0
	\end{equation}
	with 
	\begin{equation}
	\frac{d X_n}{dt} = B_n(t) X_n,
	\end{equation}
	one has
	\begin{equation}  \label{bn.1}
	B_n(t) = \e^{-\ad_{\Omega_n}} (B_{n-1} - G_n), \qquad \mbox{ with } \qquad G_n(t) = d \exp_{\Omega_n(t)} (\dot{\Omega}_n(t)) 
	\end{equation}  
	so that the solution of equation (\ref{ccv.1}) is expressed as
	\begin{equation}   \label{ntransf}
	x(t) = \e^{\Omega_1(t)} \, \e^{\Omega_2(t)} \cdots \e^{\Omega_n(t)} X_n(t).
	\end{equation}
	Alternatively, we can write $B_n(t)$ in Eq. (\ref{bn.1}) as follows. Since it is also true that \cite{blanes09tme}
	\begin{eqnarray}
	d \exp_{\Omega_n(t)} (\dot{\Omega}_n(t))  & = & \left(\frac{d }{dt} \e^{\Omega_n} \right) \, \e^{-\Omega_n} =
	\left( \int_0^1 \e^{x \Omega_n} \dot{\Omega}_n \e^{(1-x) \Omega_n} dx \right) \e^{-\Omega_n}  \nonumber \\
	& = & \int_0^1 \e^{x \, \Omega_n} \dot{\Omega}_n \e^{-x \, \Omega_n} dx, 
	\end{eqnarray}
	we have
	\begin{equation}  \label{bn.2}
	\begin{aligned}
	B_n = \ & \e^{-\Omega_n} B_{n-1} \e^{\Omega_n} - \int_0^1 \e^{-u \, \Omega_n} \dot{\Omega}_n \e^{u \, \Omega_n} du \\
	= \ & \sum_{k \ge 0} \frac{(-1)^k}{(k+1)!} \left( (k+1) \ad_{\Omega_n}^k B_{n-1} - \ad_{\Omega_n}^k \dot{\Omega}_n \right),    \qquad n \ge 1. 
	\end{aligned}
	\end{equation}  
	The important point is of course how to choose $B_n$, or alternatively $\Omega_n$, i.e., the specific requirements each transformation has to satisfy {in order to be useful to approximately solve Eq. (\ref{ccv.1})}. There are obviously many
	possibilities and in the following we analyze two of them, leading to two different and well known exponential perturbation factorizations mentioned
	in the Introduction,
	namely the Wilcox \cite{wilcox67eoa} and Fer \cite{fer58rdl} expansions.

\section{Wilcox expansion}
\label{section.3}

\subsection{Recurrences}

{Let us introduce the (dummy) parameter $\lambda$ in Eq. (\ref{ccv.1}) and   replace $A$ with $\lambda  A$.} This is helpful when collecting
coefficients, and at the end we can always take $\lambda = 1$.

Since the solution of Eq. 
(\ref{ccv.1}) when $A$ is constant, or more generally when $A(t_1) A(t_2) = A(t_2) A(t_1)$ for all $t_1 \ne t_2$, is $x(t) = \exp( \int_0^t A(u) du)$, it makes
sense to take the generator for the first transformation as
\begin{equation}
\Omega_1(t) =  \int_0^t B_0(u) \, du = \lambda \int_0^t A(u) \, du \equiv \lambda \, W_1(t).
\end{equation}   
Then, according to Eqs. (\ref{b1})-(\ref{b2}), we have
\begin{equation}
C_1 \equiv B_0 - G_1 = \sum_{k \ge 1} \lambda^k (b_{0,k} - g_{1,k}) 
\end{equation}
where
\begin{equation}
\begin{aligned} 
&  b_{0,1} = A(t), \qquad\qquad b_{0,l} = 0, \;\; l > 1 \\
&  g_{1,1} = \dot{W}_1,  \qquad\qquad \; g_{1,l} = \frac{1}{l!} \ad_{W_1}^{l-1} \dot{W}_1, \;\; l > 1.
\end{aligned}
\end{equation}   
With the choice $\dot{W}_1 = A(t)$, it turns out that $B_1$ is a power series in $\lambda$ starting with $\lambda^2$,
\begin{equation}
B_1(t) = \e^{- \ad_{\Omega_1}} (B_0 - G_1) =    \e^{-\lambda \, \ad_{W_1}} C_1 = \sum_{l \ge 2} \lambda^l \, b_{1,l},
\end{equation}
where   
\begin{equation}
b_{1,l} = \sum_{k=0}^{l-2} \frac{(-1)^k}{k!}  \ad_{W_1}^k \big( b_{0,l-k} - g_{1,l-k} \big).
\end{equation}
We can  analogously choose the second transformation proportional to $\lambda^2$, i.e., as $\Omega_2(t) \equiv \lambda^2 W_2(t)$
for a given $W_2$ to be determined.
Then a straightforward calculation shows that
\begin{equation}  \label{eq.aux.1}
C_2 \equiv B_1 - G_2 = \sum_{l \ge 2} \lambda^l (b_{1,l} - g_{2,l}) 
\end{equation}
with  
\begin{equation}
g_{2,2} = \dot{W}_2, \qquad g_{2,2l} = \frac{1}{l!} \ad_{W_2}^{l-1} \dot{W}_2, \qquad g_{2,r} =0, \quad r \ne 2l.
\end{equation}
The generator $W_2$ is then obtained by imposing that $b_{1,2} - g_{2,2} = 0$ in Eq. (\ref{eq.aux.1}), i.e.,
\begin{equation}
g_{2,2} = \dot{W}_2 =  b_{1,2} = b_{0,2} - g_{1,2} = -\frac{1}{2} \ad_{W_1} \dot{W}_1.
\end{equation}
In this way $B_2(t)$ is a power series in $\lambda$ starting with $\lambda^3$,
\begin{equation}
B_2(t) = \e^{-\lambda^2 \, \ad_{W_2}} C_2 = \sum_{l \ge 3} \lambda^l \, b_{2,l}
\end{equation}
with
\begin{equation}
b_{2,l} = \sum_{k=0}^{[(l-1)/2]-1} \frac{(-1)^k}{k!}  \ad_{W_2}^k \big( b_{1,l-2k} - g_{2,l-2k} \big),
\end{equation}
where $[\cdot]$ stands for the integer part of the argument.
In general, the $n$-th transformation $\Omega_n(t) \equiv \lambda^n \, W_n(t)$ 
is determined in such a way that the power series of $B_n(t)$ starts with $\lambda^{n+1}$. This can be done as follows: from
$B_{n-1} = \sum_{l \ge n} \lambda^l b_{n-1,l}$,
we compute
\begin{equation}  \label{algo10.1}
C_n \equiv B_{n-1} - G_n = \sum_{r = n}^{\infty} \lambda^r c_{n,r} =  \sum_{r = n}^{\infty} \lambda^r (b_{n-1,r} - g_{n,r})
\end{equation}
with
\begin{equation}  \label{algo10.2}
g_{n,r} = \left\{ \begin{array}{ll}
\displaystyle \frac{1}{l!} \ad_{W_n}^{l-1} \dot{W}_n, \qquad & r = n \, l \\
0, & r \ne n \, l 
\end{array} \right. \qquad l = 1, 2, \ldots
\end{equation}
Then, $W_n$ is obtained by taking $b_{n-1,n} - g_{n,n} = 0$, i.e., 
\begin{equation}  \label{algo10.3}
\dot{W}_n  =   b_{n-1,n}     
\end{equation}
and finally $B_n$ is determined as
\begin{equation}  \label{algo10.4}
B_n = \e^{-\lambda^n \, \ad_{W_n}} C_n = \sum_{l=n+1}^{\infty}  \lambda^l \, b_{n,l}
\end{equation}
with
\begin{equation}  \label{algo10.5}
b_{n,l} = \sum_{k=0}^{[\frac{l-1}{n}]-1} \, \frac{(-1)^k}{k!}
\ad_{W_{n}}^k c_{n,l- n \, k}.
\end{equation}
Notice that, in view of Eq. (\ref{algo10.2}) and Eq. (\ref{algo10.5}), equation (\ref{algo10.3}) simplifies to
\begin{equation}  \label{algo10.6}
\dot{W}_n  =   b_{n-2,n}    \qquad \mbox{ for } \qquad n \ge 3.     
\end{equation}
The solution of Eq. (\ref{ccv.1}) is expressed, after $n$ such transformations, as
\begin{equation}  \label{approx.W1}
x(t) = \e^{\lambda W_1(t)} \, \e^{\lambda^2 W_2(t)} \, \cdots \, \e^{\lambda^n W_n(t)} X_n(t).
\end{equation}
An approximation to the exact solution containing all the dependence up to $\mathcal{O}(\lambda^n)$ is obtained by 
taking $X_n = x_0$  in Eq. (\ref{approx.W1}).

In spite of this, the truncated factorization  $ \e^{\lambda W_1(t)} \, \e^{\lambda^2 W_2(t)} \, \cdots \, \e^{\lambda^n W_n(t)}$ still shares relevant qualitative properties, such as orthogonality, unitarity, etc.

Recursion (\ref{algo10.1})-(\ref{algo10.6}) allows one to construct \emph{any} generator $W_k$ of the expansion in terms of $W_1, \ldots, W_{k-1}$. 
{For illustration, we next collect the first terms:}
\begin{equation}  \label{wsrec}
\begin{aligned}
& \dot{W}_1  =  A(t), \qquad\quad  \dot{W}_2 =  - \frac{1}{2} \ad_{W_1} \dot{W}_1, \qquad\quad
\dot{W}_3  =   \frac{1}{3} \ad_{W_1}^2  \dot{W}_1 \\
&  \dot{W}_4  =    -\frac{1}{8} \ad_{W_1}^3  \dot{W}_1 - \frac{1}{2} \ad_{W_2}  \dot{W}_2, \qquad\quad 
\dot{W}_5  =  \frac{1}{30} \ad_{W_1}^4 \dot{W}_1 - \ad_{W_2} \dot{W}_3
\end{aligned}
\end{equation}
{This is the way the Wilcox expansion is built up.}

\subsection{Explicit expressions for $W_n(t)$}
\label{sec.explicit}

Although the recursive procedure (\ref{algo10.1})-(\ref{algo10.6}) turns out to be very computationally efficient to construct  the exponents $W_n(t)$ for a given $A(t)$ in practice, it is clear that much insight about the expansion can be gained if an explicit expression for \emph{any} $W_n$ can be constructed, thus generalizing the treatment done originally by
Wilcox up to $n=4$ \cite{wilcox67eoa}.

Such an expression could be obtained, in principle, by working out the recurrence (\ref{algo10.1})-(\ref{algo10.6}), but 
a more direct approach consists in
comparing the Dyson perturbation series of $U(t)$ \cite{galindo90qme} in the associated initial value problem
\begin{equation}
\dot{U} = \lambda \, A(t) U,   \qquad U(0) = I,
\end{equation}
i.e.,    
\begin{equation}
U(t) = I + \sum_{n=1}^{\infty} \lambda^k P_k(t), \qquad \mbox{ with } \qquad P_k(t) = \int_0^t dt_1 \cdots \int_0^{t_{k-1}} dt_k \, A(t_1) \cdots A(t_k)
\end{equation}
with the expansion in $\lambda$ of the factorization
\begin{equation}
U(t) = \e^{\lambda W_1(t)} \, \e^{\lambda^2 W_2(t)} \, \cdots \, \e^{\lambda^n W_n(t)}  \cdots.
\end{equation}
Thus, for the first terms one has
\begin{equation}
\begin{aligned}
P_1 & =W_1, \qquad\quad   P_2  = W_2 + \dfrac 12 W_1^2, \qquad\quad  
P_3  =W_3 + W_1W_2+\frac{1}{3!} W_1^3,\\
P_4 & =W_4 + W_1 W_3+ \frac{1}{2} W_1^2 W_2 + \dfrac 12 W_2^2 + \frac{1}{4!} W_1^4.
\end{aligned}
\end{equation}
In general, we can write
\begin{equation}  \label{p.w}
P_n = \sum_{p(n) \atop  i_1 \le i_2 \le \cdots \le i_k}  \frac{1}{r_k!} \, W_{i_1} W_{i_2} \cdots W_{i_k}, \qquad i_1 \le i_2 \le \cdots \le i_k,
\end{equation}   
where the sum is extended over the total number of partitions $p(n)$ of the integer $n$. {We recall that a partition $p(n)$ of the integer $n$
	is an $n$-tuple $(i_1, i_2, \ldots, i_k)$ such that $i_1 + i_2 + \cdots + i_k = n$, with ordering $i_1 \le i_2 \le \cdots \le i_k$. Thus,  the
	seven partitions of $n=5$ with the chosen ordering are $(5)$, $(1,4)$, $(2,3)$, $(1,1,3)$, $(1,2,2)$, $(1,1,1,2)$ and $(1,1,1,1,1)$.
	In Eq. (\ref{p.w}), 
	$r_k$ is the number of repeated indices in the partition considered.}

By working out equation (\ref{p.w}), one can invert the relations and express $W_n$ in terms of $P_1, \ldots, P_{n-1}$ {for any $n \ge 1$}. Thus, i  one obtains
\begin{equation}
\begin{aligned}
W_1 & =P_1, \qquad\quad  W_2   = P_2-\frac 12 P_1^2, \qquad\quad  
W_3  = P_3- P_1 P_2+\frac{1}{3} P_1^3, \\
W_4  & = P_4 - P_1 P_3 + \frac 34 P_1^2 P_2+\frac 14 P_2 P_1^2-\frac 12 P_2^2-\frac 14 P_1^4.
\end{aligned}
\end{equation}
Notice that $W_n$, $n \ge 2$, is expressed in terms of products of iterated integrals $P_{i_1} \cdots P_{i_j}$. {Interestingly, it is  possible to express these products as proper time-ordered integrals by
	using a procedure developed  in \cite{arnal18agf}.} If we denote
\begin{equation}   \label{permu.it}
A(i_1 i_2 \ldots i_n ) \equiv \int_0^t dt_1 \int_0^{t_1} dt_2 \cdots \int_0^{t_{n-1}} dt_n \, A(t_{i_1}) A(t_{i_2}) \cdots A(t_{t_n}),
\end{equation}
so that $P_n(t) = A(1 2 \ldots n)$, then 
\begin{equation}
\begin{aligned}
& W_2 = A( 1 2 ) - \frac{1}{2}  A( 1 ) \cdot  A( 1 ) \\
& W_3 = A(1 2 3) - A(1) \cdot A(1 2)  + \frac{1}{3} A( 1 )  \cdot A( 1 ) \cdot A(1),
\end{aligned}
\end{equation}
etc. Taking into account Fubini's theorem, 
\begin{equation}\label{intxy}
\int_0^\alpha dy \int_y^\alpha f(x,y) \, dx =
\int_0^\alpha dx \int_0^x f(x,y) \, dy,  
\end{equation}
 it is clear that $A( 1 ) \cdot  A( 1 )  =  A( 1 2 )  +  A( 2 1 )$, and thus
\begin{equation}  \label{w2}
W_2 =  A( 1 2 ) - \frac{1}{2} \big( A( 1 2 ) + A( 2 1 ) \big) = \frac{1}{2} \big(  A(1 2 ) -  A( 2 1 ) \big). 
\end{equation}
 We can proceed analogously with the following products
\begin{equation}  \label{prod.1}
\begin{aligned}
&    A( 1 ) \cdot A(1 2 ) = A(1 2 3 ) + A( 2 1 3 )  + A( 3 1 2 )  \\
& A( 1 ) \cdot A( 1 ) \cdot A( 1 ) =  A( 1 2 3 ) + A( 1 3  2 )  + A( 2 1 3 ) +
A( 2 3 1 )  + A( 3 1 2 )  + A( 3 2 1 ), 
\end{aligned}
\end{equation}
 so that
\begin{equation}  \label{w3}
W_3 =  \frac{1}{3} \, A(1 2 3) + \frac{1}{3} \, A(1 3 2)  - \frac{2}{3} \, A(2 1 3) + \frac{1}{3} \, A(2 3 1)  - \frac{2}{3}  \, A(3 1 2) + 
\frac{1}{3} \, A(3 2 1).  
\end{equation}

Carrying out this argument to any order, we can expand all the  products of integrals appearing in $W_n$. As a result, each product is
	replaced by the sum of all possible permutations of time ordering consistent with the time ordering in the factors of this product \cite{dragt83con}. 

At this point, it is illustrative to consider in detail some  examples. Thus, the product $A( 1 ) \cdot A( 1 2 )$ gives the sum of 
	all permutations of three elements such 
	that the second index is less than the third
	one. With respect to $A( 1 ) \cdot A( 1 ) \cdot A( 1 )$, since there is no special ordering, then all possible permutations have to be taken into account. 
Finally, for the product $P_1 P_3$ appearing in $W_4$ one has
\begin{equation}
P_1 P_3 = A( 1 ) \cdot A( 1 2 3 )  = A( 4 1 2 3) + A( 3 1 2 4 ) + A( 2 1 3 4 )  + A( 1 2 3 4 ).
\end{equation}
Proceeding in a similar way, one can show that any product of iterated integrals can be expressed as a sum of iterated integrals. This property
	is in fact related with a much deeper characterization of the group of permutations \cite{arnal18agf}. If $\mathfrak{S}Symm$ denotes the graded
$\mathbb{Q}$-vector space with fundamental basis given by the disjoint union of the symmetric groups $S_n$ for all $n \ge 0$, then
it is possible to define a product $\ast$ of permutations and a coproduct 
$\delta$ in $\mathfrak{S}Symm$ so that there exists a 
one-to-one correspondence between iterated integrals and permutations:
\begin{equation}
A(\sigma) \cdot A(\tau) = A(\sigma \ast \tau).
\end{equation}  
The product $\ast$ was introduced in \cite{agrachev94tsp}, 
and together with the coproduct $\delta$ endow $\mathfrak{S}Symm$ with a structure of Hopf algebra \cite{hazewinkel10ara},
the so-called Malvenuto--Reutenauer Hopf algebra of permutations \cite{malvenuto95dbq}.

In sum, the \emph{general} structure of the Wilcox expansion terms reads 
\begin{equation}\label{eq:Wgeneral}
W_n=\sum_{k=1}^{n!} \omega_k^{(n)}A(\sigma_k),
\end{equation} 
where the summation extends over \textit{all} the $n!$ permutations $\sigma_k \in S_n$ of $\{1,2, \ldots, n\}$.  
The weights $\omega_k^{(n)}$ are  given by rational numbers {that can be determined algorithmically for any $n$, although, it seems not obvious
	to find a general expression for them, in contrast with the Magnus expansion, for which such a closed formula exists \cite{strichartz87tcb}. 
	This can be then
	considered as an open problem.} 

Moreover, if one
is interested in getting a compact expression for $W_n$ in terms of independent nested commutators of $A(t)$, as is done in \cite{wilcox67eoa} up
to $n=4$, one can use the class of bases proposed by Dragt \& Forest in \cite{dragt83con} for the Lie algebra
generated by the operators $A(t_1), A(t_2), \ldots, A(t_n)$.  The same procedure carried out in \cite{arnal18agf} for the Magnus expansion can be
applied here, so that one gets the general formula
\begin{equation}   \label{wil.com3}
\begin{aligned}
& W_n(t) =  \sum_{\tau_k}  c_{\tau_k}^{(n)}
\int_0^t dt_1 \int_0^{t_1} dt_2 \cdots \int_0^{t_{n-1}} dt_n \,  \\
&  \qquad\qquad\qquad  [A(t_{\tau(2)}), [A(t_{\tau(3)}) \cdots 
[A(t_{\tau(n)}), A(t_1)] \cdots ]].
\end{aligned}    
\end{equation}
Here the sum extends over the $(n-1)!$ permutations $\tau_k$ of the elements $\{2, 3, \ldots, n \}$ and $c_{\tau_k}^{(n)}$ is 
a rational number that depends on the
particular permutation. For a given permutation, say $\tau_k$, its value coincides with the prefactor in Eq. (\ref{eq:Wgeneral}) of the particular
term $A(\sigma_k)$ corresponding to the permutation $\sigma_k \in S_n$ such that
\begin{equation}
\{\sigma_k(1), \sigma_k(2), \ldots, \sigma_k(n-1),\sigma_k(n) \} = \{ \tau_k(2), \tau_k(3), \ldots, \tau_k(n), 1 \}.
\end{equation}   

Thus, if we denote 
\begin{equation}
\label{formfs}
A[i_1 i_2 \ldots i_n] \equiv \int_0^t dt_1 \int_0^{t_1} dt_2 \cdots \int_0^{t_{n-1}} dt_n [A(t_{i_1}), [A(t_{i_2}), \cdots [A(t_{i_{n-1}}), A(t_{i_n})] \cdots ]]
\end{equation}
we get {for the first terms
	\begin{equation}
	\begin{aligned}
	& W_2 = -\frac{1}{2} A[21], \\
	&  W_3 =  \frac{1}{3} A[231]+\frac{1}{3} A[321], \\
	& W_4  = -\frac{1}{4}A[3241]-\frac{1}{4}A[4231]-\frac{1}{4}A[4321], \\
	& W_5 = -\frac{2}{15} A[23451] - \frac{2}{15} A[23541] - \frac{2}{15} A[24351] - \frac{2}{15} A[24531] \\
	& \qquad - \frac{2}{15} A[25341] - \frac{2}{15} A[25431] + \frac{1}{5} A[32451] + \frac{1}{5} A[32541] - \frac{2}{15} A[34251] \\
	& \qquad - \frac{2}{15} A[34521] - \frac{2}{15} A[35241] - \frac{2}{15} A[35421] + \frac{1}{5} A[42351] + \frac{1}{5} A[42531] \\
	& \qquad + \frac{1}{5} A[43251] + \frac{1}{5} A[43521] - \frac{2}{15} A[45231] - \frac{2}{15} A[45321] + \frac{1}{5} A[52341] \\
	& \qquad + \frac{1}{5} A[52431] + \frac{1}{5} A[53241] + \frac{1}{5} A[53421] + \frac{1}{5} A[54231] + \frac{1}{5} A[54321] 
	\end{aligned}
	\end{equation}
}
In sum, the general structure of the Wilcox expansion terms \textit{via} commutators 
uses the same weights as in Eq. (\ref{eq:Wgeneral}) and reads 
\begin{equation}\label{eq:Wcommut}
W_n=\sum_{k=1}^{(n-1)!}{\vphantom{\sum}}' \omega_k^{(n)}A[\sigma_k],
\end{equation}
where the primed sum  requires the rightmost element in the permutation $\sigma(k)$ to be invariant (as in Eq. (\ref{wil.com3})). This element may be chosen at will  and, whatever it be that value, the permutations are build up  with the remaining $n-1$ elements. Different, but equivalent, expressions for $W_n$ in terms of commutators are obtained  depending on the value fixed at the rightmost position. {We stress once again that, although only the first
	terms have been collected here for simplicity, the whole procedure is algorithmic in nature and has been implemented in a computer algebra system
	furnishing to evaluate explicitly $W_n(t)$ for any $n$ \cite{wilcoxcode}. Notice that $W_n(t)$ involves a linear combination of $(n-1)!$ iterated integrals.}

\subsection{{Convergence of Wilcox expansion}}
\label{conv.expansion}

Recursion (\ref{algo10.1})-(\ref{algo10.6}) is also very useful to provide estimates for the radius of convergence of the Wilcox expansion when 
Eq. (\ref{ccv.1}) is defined in a 
Banach algebra $\mathcal{A}$, i.e.,
an algebra that is also a complete normed linear space with a sub-multiplicative norm,
\begin{equation}
\|X \, Y \| \le \|X\| \, \|Y\|.
\end{equation}
If this is the case, then $\| \ad_X Y  \| \le 2 \, \|X\| \, \|Y\|$ and, in general, $\| \ad_X^n Y\| \le 2^n \|X\|^n \, \|Y\|$.

As shown in \cite{bayen79otc,arnal17ots}, if the series
\begin{equation}
M(\lambda;t) = \sum_{j=1}^{\infty} \lambda^j \, \|W_j(t)\|
\end{equation}
has a certain radius of convergence $r_c$ for a given $t$, then, for $\lambda  < r < r_c$, the sequence of functions
\begin{equation}
\Psi_n \equiv \e^{\lambda W_1(t)} \,  \e^{\lambda^2 W_2(t)} \, \cdots \, \e^{\lambda^n W_n(t)}
\end{equation}
converges uniformly on any compact  subset of the ball $B(0,r_c)$. Thus, studying the convergence of the Wilcox expansion reduces to
analyze the series $M(\lambda;t)$ and in particular its radius of convergence $r_c$. 

Let $k(t)$ be a function such that $\|A(t)\| \le k(t)$ and denote $K(t) = \int_0^t k(s) \, ds$. Then, clearly
\begin{equation}
\|b_{0,1}\| \le k(t); \qquad\quad \|b_{0,l}\| = 0, \quad l=2,3,\ldots
\end{equation}      
and
\begin{equation}
\|W_1\|   \le   K(t), \qquad \|W_2\| \le \frac{1}{2} K^2(t).
\end{equation}
In general,  the following bounds can be established by induction:
\begin{eqnarray}
\|g_{n,r}(t) \| & \le & \beta_{n,r} \, K^{r-1}(t) \, k(t), \nonumber \\
\|b_{n,l}(t)\| & \le & \alpha_{n,l} \, K^{l-1}(t) \,  k(t), \qquad\qquad n=1,2,\ldots; \quad l > n \nonumber \\
\|W_n(t)\| & \le &  c_n \, K^n(t),
\end{eqnarray}
where
\begin{eqnarray}
\alpha_{0,1} & = & 1; \qquad \alpha_{0,l} = 0, \;\; l > 1 \nonumber \\
\alpha_{n,l} & = &    \sum_{j=0}^{[\frac{l-1}{n}]-1} \frac{1}{j!} \, 2^j \, c_n^j  \, (\alpha_{n-1,l-nj} + \beta_{n,l-nj}) \quad\qquad n \ge 1, \; l > n \nonumber \\
\beta_{n,r}  & =  & \left\{  \begin{array}{ll}
\displaystyle \frac{1}{l!} \, 2^{l-1} \, n \, c_n^l, \quad & \left[ \frac{r}{n} \right] = \frac{r}{n} = l \\
0,  & \left[ \frac{r}{n} \right] \ne \frac{r}{n}
\end{array} \right. \\
c_1 & = & 1, \quad c_2 = \frac{1}{2}, \qquad 
c_n  =  \frac{1}{n} \, \alpha_{n-2,n},  \qquad n > 2. \nonumber
\end{eqnarray}

\begin{figure}[ht]
	\begin{center}
		\includegraphics[scale=0.45]{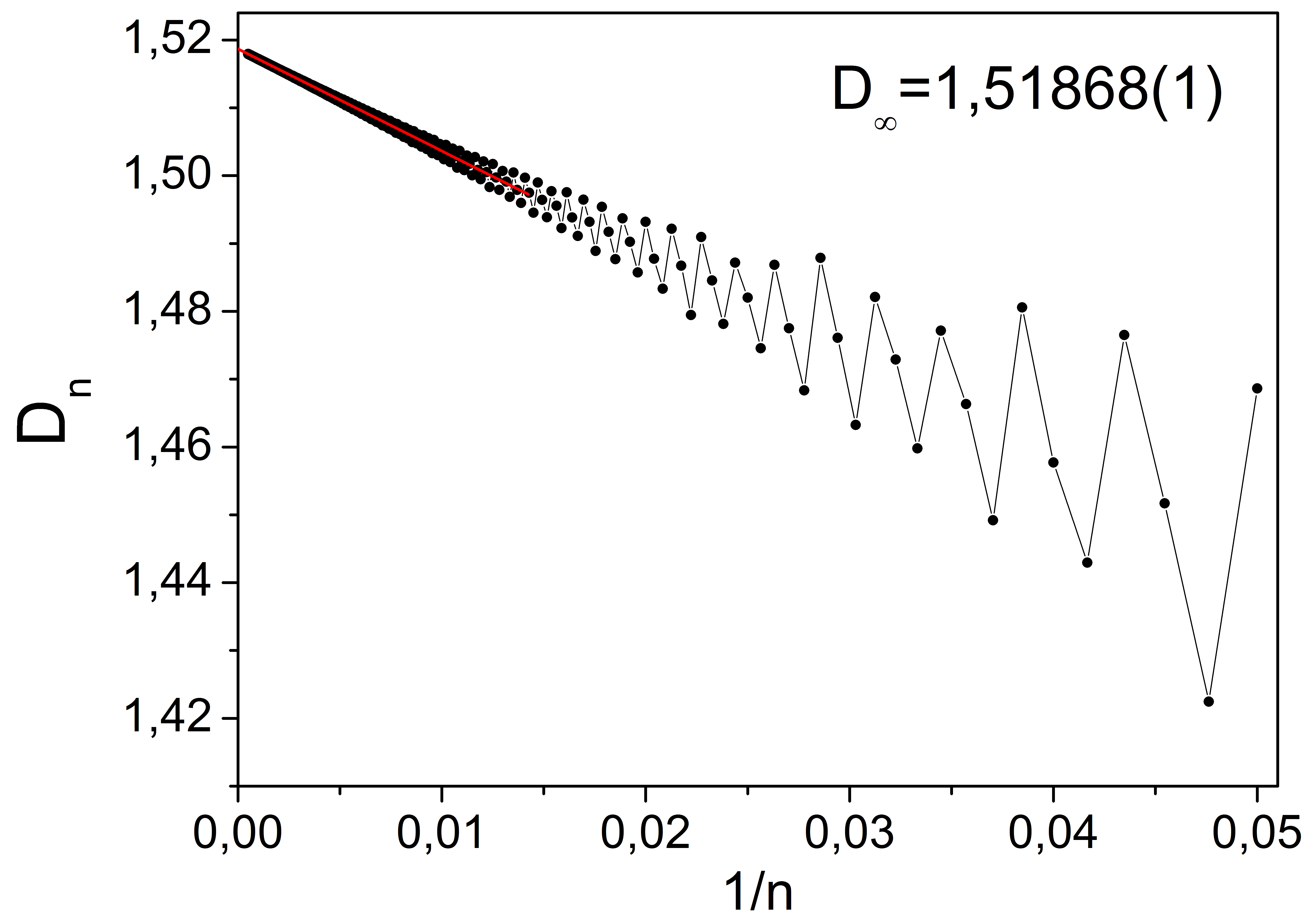}
		\caption{{\small $D_n$ as a function of $1/n$, and linear extrapolation (red line).}} \label{figu.conv}
	\end{center}
\end{figure}

It is clear that if the series $\sum_{j \ge 1} c_j \, K^j(t)$ converges, so does $M(\lambda=1;t)$. 
Therefore, a sufficient condition for convergence of the Wilcox
expansion is obtained by imposing
\begin{equation}
\lim_{n \rightarrow \infty} \frac{c_{n+1} K^{n+1}(t)}{c_n K^n(t)}  =  K(t) \lim_{n \rightarrow \infty} D_n < 1,
\end{equation}   
where
\begin{equation}
D_n \equiv \frac{n}{n+1} \frac{\alpha_{n-1,n+1}}{\alpha_{n-2,n}}.
\end{equation}
We have computed this quantity up to $n=2000$  and then extrapolated to the limit $(1/n) \rightarrow 0$. Then $D_n \rightarrow D_{\infty} = 1.51868$,
as seen in Figure \ref{figu.conv},
and thus the convergence of the Wilcox expansion is ensured at least for values of time $t$ such that
\begin{equation}  \label{conv.wil}
\int_0^t \|A(s)\| ds \le K(t) < \xi_W = \frac{1}{D_{\infty}} \approx 0.65846
\end{equation}   
This type of extrapolation has {also} been used to  estimate the convergence radius of the Magnus expansion \cite{moan01cot}.
{Although the estimate is not completely analytic, the same type of computation has provided accurate results in other settings. In particular,
	for the Magnus expansion such an estimate fully agrees with a purely theoretically deduced bound \cite{blanes98maf,blanes09tme,moan01cot}.}

\section{Fer{--like} expansions}

\subsection{Standard Fer expansion}

{In} forming the Wilcox expansion the first transformation is chosen in such a way that $\dot{\Omega}_1 =  B_0(t) $, whereas 
$\dot{\Omega}_n \ne B_{n-1}$ for $n \ge 2$. It makes sense, then, to analyze what happens if we impose this condition at each step of the procedure
\begin{equation}
\dot{\Omega}_n(t) = B_{n-1}(t) \qquad \mbox{ or, equivalently, } \qquad \Omega_n(t) = \int_{0}^t B_{n-1}(u) du
\end{equation}
for all $n \ge 1$. 
In that way, expression (\ref{bn.2}) for $B_n$ clearly simplifies to
\begin{equation}   \label{bnfer}
B_n(t) = \sum_{k \ge 1} \frac{(-1)^k k}{(k+1)!} \, \ad_{\Omega_n}^k B_{n-1}, \qquad n \ge 1.
\end{equation} 
In doing so we recover precisely the Fer expansion, see \cite{fer58rdl,blanes98maf}. Again, after $n$ transformations, we get
\begin{equation}  \label{Fer}
x(t) =  \e^{\Omega_1(t)} \, \e^{\Omega_2(t)} \cdots \e^{\Omega_n(t)} X_n(t)
\end{equation}
so that if we impose $X_n(t) = x_0$ we are left with another approximation to the exact solution. Notice that 
this approximation clearly
differs from the previous Wilcox expansion for $n \ge 2$, as can be seen 
by analyzing the dependence on $\lambda$ of each transformation. Whereas $\Omega_n = \lambda^n W_n(t)$ for
the Wilcox expansion, now $\Omega_n$ contains terms of order $\lambda^{2^{n-1}}$ and higher. This can be easily shown by induction: $\Omega_1$ is
proportional to $ \lambda$, so that $B_1$, according to Eq. (\ref{bnfer}) contains terms of order $\lambda^{2}$ (coming from $[\Omega_1, B_0]$) and higher. 
In general, $B_{n-2}$ and $\Omega_{n-1}$ contain terms of order $\lambda^{2^{n-2}}$ and higher,   so the first term in the series (\ref{bnfer}) for $B_{n-1}$,
i.e., the commutator $[\Omega_{n-1}, B_{n-2}]$, produces a term of order $(\lambda^2)^{2^{n-2}} = \lambda^{2^{n-1}}$ in $\Omega_n$.

Alternatively, expressing Eq. (\ref{bnfer})  as
\begin{equation}
B_n(t) = \int_0^1 dx \int_0^x du \, \e^{-(1-u) \Omega_n} \, [ B_{n-1}, \Omega_n] \, \e^{(1-u) \Omega_n}
\end{equation}
and taking norms, it is then possible to show that the Fer expansion converges for values of $t$ such that \cite{blanes98maf}
\begin{equation}
\int_0^t \|A(s)\| ds < 0.8604065.
\end{equation}

\subsection{Intermediate Fer-like expansions}

Notice that the $\lambda$-power series of $\Omega_n$ in the Fer expansion contains infinite terms starting with
$\lambda^{2^{n-1}}$, but the corresponding truncated factorization obtained from Eq. (\ref{Fer}) by taking $X_n(t) = x_0$
is correct only up to terms of order $\lambda^{2^{n-1}}$. One might then consider yet another sequence of transformations so that each $\Omega_k$ contains
only those relevant terms leading to a correct approximation up to this order. Of course, both factorizations would be different, but nevertheless they would produce
the correct power series up to order $\lambda^{2^{n-1}}$. The corresponding factorization can be properly called a modified Fer expansion.

Our starting point is once again equation (\ref{bn.2}). Clearly, the first transformation is the same as in Fer (and Wilcox), i.e., 
\begin{equation}  \label{mf1}
\dot{\Omega}_1 = B_0 = \lambda A(t),
\end{equation}   
and thus 
\begin{equation}   \label{eq.b1}
B_1 = \sum_{k=0}^{\infty} \frac{(-1)^k k}{(k+1)!} \ad_{\Omega_{1}}^k B_0 =  \mathcal{O}(\lambda^2) ,
\end{equation}
{where the rightmost term points out the lowest $\lambda$ contribution in the sum.}

Next, to reproduce the same dependence on $\lambda$ as the Fer expansion, we need to enforce that $B_2 = \mathcal{O}(\lambda^4)$, and the question is how to
choose $\Omega_2$ guaranteeing this feature. An analysis of Eq. (\ref{bn.2}) with $n=2$ reveals that 
this is achieved by taking $\dot{\Omega}_2$ as the sum of terms in $B_1$ in Eq. (\ref{eq.b1}) contributing to $\lambda^2$ and $\lambda^3$, i.e.,
\begin{equation}  \label{mf2}
\dot{\Omega}_2 = -\frac{1}{2} \ad_{\Omega_1} B_0 + \frac{1}{3} \ad_{\Omega_1}^2 B_0,
\end{equation}
since the next term appearing in the expression of $B_2$ involves the computation of  $\ad_{\Omega_1}^3 B_0  = \mathcal{O}(\lambda^4)$. Thus
\begin{equation}  \label{b2mf}
B_2 = \sum_{k=1}^{\infty} \frac{(-1)^k}{(k+1)!} \left( (k+1) \ad_{\Omega_2}^k B_1 - \ad_{\Omega_2}^k \dot{\Omega}_2 \right) = \mathcal{O}(\lambda^4).
\end{equation}
Likewise, $\Omega_3$
is to be designed so that 
\begin{equation}
B_3 = \sum_{k=0}^{\infty} \frac{(-1)^k}{(k+1)!} \left( (k+1) \ad_{\Omega_3}^k B_2 - \ad_{\Omega_3}^k \dot{\Omega}_3 \right) = 
\mathcal{O}(\lambda^8),
\end{equation}
and this is guaranteed by taking $\Omega_3$ as the sum all the terms in $B_2$ contributing to powers from $\lambda^4$ up to
$\lambda^7$. From Eq. (\ref{b2mf}) it is clear that
\begin{equation}
B_2 = \underbrace{-\frac{1}{2} \left( 2 \ad_{\Omega_2} B_1 - \ad_{\Omega_2} \dot{\Omega}_2 \right)}_{\mathcal{O}(\lambda^4)} + 
\underbrace{\frac{1}{3!} \left( 3 \ad_{\Omega_2}^2 B_1 -
	\ad_{\Omega_2}^2 \dot{\Omega}_2 \right)}_{\mathcal{O}(\lambda^6)} + \mathcal{O}(\lambda^8),
\end{equation}
where only the relevant terms in the expansion in $B_1$ have to be taken into account. In this way, we can take
\begin{equation}  \label{mf3}
\dot{\Omega}_3 = - \ad_{\Omega_2} B_1^{[4]} + \frac{1}{2} \ad_{\Omega_2} \dot{\Omega}_2 + \frac{1}{2} \ad_{\Omega_2}^2 B_1^{[2]} - \frac{1}{6} 
\ad_{\Omega_2}^2 \dot{\Omega}_2
\end{equation}
with
\begin{equation}
B_1^{[j]} \equiv \sum_{k=1}^j \frac{(-1)^k}{(k+1)!} k \, \ad_{\Omega_1}^k B_0.
\end{equation}   
Notice, however, that since the second term in $\Omega_2$ in Eq. (\ref{mf2}) is $\mathcal{O}(\lambda^3)$, the expression (\ref{mf3}) does contain some
contributions in $\lambda^8$ and $\lambda^9$ that in principle could be removed. We prefer, however, to maintain them just to have a more compact
expression.

For this modified Fer expansion, $\Omega_n$ is chosen in general
so that $\dot{\Omega}_n$ is precisely the sum of all terms of $B_{n-1}$ containing terms of powers from
$\lambda^{2^{n-1}}$ up to  $\lambda^{2^n-1}$ and then truncating appropriately the series of $B_{n-2}, \ldots, B_1$.      

Other possibilities for choosing $B_k$ at the successive stages clearly exist, and according with the particular election, different intermediate
Fer-like expansions result. {In practice, one of those combinations of commutators could be more easily computed for a specific problem.}

\section{Applications}

\subsection{{Wilcox expansion as the continuous analogue of Zassenhaus formula}}

The Zassenhaus formula may be considered as the dual of the Baker--Campbell--Hausdorff (BCH) formula \cite{casas09aea} 
in the sense that it relates the
exponential of the sum of two non-commuting operators $X$ and $Y$ with an infinite product of exponentials of these operators and their nested commutators. More
specifically, 
\begin{equation}  \label{zass.1}
\e^{X + Y} = \e^X \, \e^Y \, \prod_{n=2}^{\infty} \e^{C_n(X,Y)} = \e^X \, \e^Y \, \e^{C_2(X,Y)} \,
\e^{C_3(X,Y)} \, \cdots  \, \e^{C_k(X,Y)} \, \cdots,
\end{equation}
where $C_k(X,Y)$ is a homogeneous Lie polynomial in $X$ and $Y$ of degree $k$
\cite{casas12eco,magnus54ote,suzuki77otc,weyrauch09ctb,wilcox67eoa}. A very efficient procedure to generate all the terms in Eq. (\ref{zass.1}) 
is presented in \cite{casas12eco} and allows to construct $C_n$ up to a prescribed value of $n$ directly in terms of the minimum number of 
independent commutators involving $n$ operators $X$ and $Y$.

In view of the formal similarity between Eq. (\ref{approx.W1}) and Eq. (\ref{zass.1}), Wilcox expansion has been also described as the 
``continuous analogue of the Zassenhaus formula" \cite{blanes98maf}, just as the Magnus expansion is sometimes called
the continuous version of the BCH formula.  To substantiate this claim, we next reproduce the Zassenhaus formula (\ref{zass.1})
by applying the procedure of Section
\ref{section.3} to a particular initial value problem, namely the abstract equation
\begin{equation}  \label{zass.2}
\dot{U} = \lambda (X + Y) \, U, \qquad U(0) = I,
\end{equation}
where $X$ and $Y$ are two non-commuting constant operators. 

The formal solution is of course
$U(t) = \e^{t \lambda (X+Y)}$, but we can also solve Eq. (\ref{zass.2}) by first integrating $\dot{U}_0 = \lambda X U_0$ and factorizing $U(t)$ as
$U(t) = U_0 \, U_I = \e^{t \lambda X} \, U_I$, where $U_I$ obeys the equation
\begin{equation}  \label{zass.3}
\dot{U}_I = \lambda \, \e^{-t \lambda X} Y \e^{t \lambda X} \, U_I \equiv A_{\lambda}(t) \, U_I,
\end{equation}   
and finally apply to Eq. (\ref{zass.3}) the sequence of transformations leading to the Wilcox expansion. Notice, however, that now the coefficient matrix
$A_{\lambda}(t)$ is an infinite series in $\lambda$,
\begin{equation}
A_{\lambda}(t) = \lambda \, \e^{- t \lambda \ad_{X}} Y = \sum_{j\ge 0} \frac{(-1)^j}{j!} t^j \lambda^{j+1} \ad_X^j Y,
\end{equation}
so that, when applying the recursion (\ref{algo10.1})-(\ref{algo10.6}), $\dot{\Omega}_1$ is no longer $B_0 \equiv  A_{\lambda}(t)$, but the term 
in $A_{\lambda}(t)$ which is proportional to $\lambda$.
In other words,
\begin{equation}
\dot{W}_1 =  Y, \qquad \mbox{ and thus } \qquad W_1(t) = t  \, Y.
\end{equation}
{After some computation, one arrives at}
\begin{equation}  \label{b0l}
b_{0,l} = \frac{(-1)^{l-1}}{(l-1)!} t^{l-1} \ad_X^{l-1} Y.
\end{equation}
Since $\dot{W}_1 = b_{0,1} = Y$, then clearly $g_{1,l} = 0$ for all $l > 1$ and
\begin{equation}
b_{1,l} = \sum_{k=0}^{l-2} \frac{(-1)^k}{k!}  \ad_{W_1}^k  b_{0,l-k}. 
\end{equation}
By imposing  $b_{1,2} = g_{2,2} =\dot{W}_2$, we get
\begin{equation}
\dot{W}_2 =  -t \, \ad_X Y, \qquad \mbox{ so that } \qquad W_2(t) = -\frac{1}{2} \,  t^2 \, \ad_X Y.
\end{equation}
In general, $g_{n,n} = \dot{W}_n$, $g_{n,r} = 0$ when $r \ne n$ and the recurrence (\ref{algo10.1})-(\ref{algo10.6}) reads now
\begin{eqnarray}   \label{zass.rec}
\dot{W}_1 & = & b_{0,1} \nonumber \\
b_{n,l} & = & \sum_{k=0}^{[\frac{l-1}{n}]-1} \, \frac{(-1)^k}{k!}
\ad_{W_{n}}^k  b_{n-1,l-nk}, \qquad n=1,2,\ldots \\
\dot{W}_n & = &  b_{n-2,n} \qquad n \ge 2  \nonumber  
\end{eqnarray}
together with Eq. (\ref{b0l}). Working out this recursion we obtain, for the first terms
	\begin{eqnarray}
	W_3(t) & = & \frac{1}{6} t^3 \ad_X^2 Y + \frac{1}{3} t^3 \ad_Y \ad_X Y  \nonumber \\
	W_4(t) & = & -\frac{1}{24} t^4 \ad_X^3 Y  - \frac{1}{8} t^4 \ad_Y \ad_X^2 Y  - \frac{1}{8} t^4 \ad_Y^2 \ad_X Y \nonumber \\
	W_5(t) & = & \frac{1}{120} t^5 \ad_X^4 Y+ \frac{1}{30} t^5 \ad_Y \ad_X^3 Y + \frac{1}{20} t^5 \ad_Y^2 \ad_X^2 Y \\
	& & + \frac{1}{30} t^5 \ad_Y^3 \ad_X Y + \frac{1}{20} t^5 \ad_{[X,Y]} \ad_X^2 Y + \frac{1}{10} t^5 \ad_{[X,Y]} \ad_Y \ad_X Y \nonumber
	\end{eqnarray}  

{One can see that} this procedure {agrees with the} algorithm {presented} in \cite{casas12eco} for {every} term $W_n$, $n \ge 1$,  in
\begin{equation}
U(t) = \e^{t \lambda (X+Y)} = \e^{t \lambda X} \, \e^{\lambda W_1(t)} \, \e^{\lambda^2 W_2(t)} \, \e^{\lambda^3 W_3(t)} \, \cdots.
\end{equation}
{The} Zassenhaus formula is recovered by taking $t=1$, {i.e., $C_n(X,Y) = W_n(t=1)$.}

\subsection{Expanding the exponential $\exp(A + \varepsilon B)$}

Bellman, in his classic book \cite{bellman72pti},  states that ``one of the great challenges of modern physics is that of obtaining useful
approximate relations  for $\e^{(A + \varepsilon B)t}$ in the case where $A B \ne B A$''. One such approximation was   {proposed and left undisclosed} in \cite[p. 175]{bellman70itm}.   Assuming that  $\e^{A + \varepsilon B}$ can be written in the form
\begin{equation}   \label{bell.0}
\e^{A + \varepsilon B} = \e^A \, \e^{\varepsilon C_1} \, \e^{\varepsilon^2 C_2} \, \e^{\varepsilon^3 C_3} \cdots ,
\end{equation}
{Bellman proposed to determine the first three terms $C_1$, $C_2$, $C_3$, and pointed out}  that,  contrary to other expansions, the product expansion
(\ref{bell.0}) is
unitary if $A$ and $B$ are skew-Hermitian.

It turns out that the Wilcox expansion can be used to provide explicit expressions for $C_n$  for any two indeterminates $A$ and $B$, as we will
see in the sequel.

{Before proceeding, it is important to remark that this problem differs from the Zassenhauss formula, in the sense that the expansion parameter
	affects only one of the operators in the exponential. The solution goes as follows. We write}
\begin{equation}  \label{bell1}
U(t) \equiv \e^{t (A + \varepsilon B)} = \e^{t A} \, V,
\end{equation}
and solve the differential equation satisfied by $V$
\begin{equation}  \label{bell2}
\frac{dV}{dt} = \varepsilon e^{-t A} B e^{t A} \, V \equiv \varepsilon \tilde{B}(t) V, \qquad V(0) = I
\end{equation}
with the Wilcox expansion, so that 
\begin{equation}   \label{bell.2a}
V(t) = \e^{\varepsilon W_1(t)} \, \e^{\varepsilon^2 W_2(t)} \, \e^{\varepsilon^3 W_3(t)} \cdots.
\end{equation}
The operators $C_i$ in Eq. (\ref{bell.0}) are then obtained by taking $t=1$.

The successive terms $W_j(t)$ in Eq. (\ref{bell.2a}) can be determined either by the recursion (\ref{algo10.1})-(\ref{algo10.6}) or
the explicit expression (\ref{wil.com3}). For the first term we get
\begin{equation}  \label{bell.2b}
W_1(t)  =  \int_{0}^{t}\tilde{B}(s) dt  = \int_0^t \e^{-s \, \ad_A} B \, ds =
\sum _{k=0}^\infty \frac{(-1)^k \,t^{k+1} }{(k+1)!} \ad_A^k B
\end{equation}
and so the expression for  $C_1$ is given by
\begin{equation}  \label{bell.2c}
C_1 = \frac{1- \e^{ \ad_A}}{\ad_A} B   = B - \frac{1}{2} [A,B] + \frac{1}{3!} [A,[A,B]] - \frac{1}{4!} [A,[A,[A,B]]] + \cdots.
\end{equation}
Although it is possible in principle to construct {explicit expressions} for $W_2$, $W_3$, etc., it is perhaps more convenient to apply the recursion
(\ref{algo10.1})-(\ref{algo10.6}) for each particular application.

\subsection{Illustrative examples}
We next particularize the \textit{Bellman problem} (\ref{bell.0}) to matrices  where
closed expressions for $C_1, C_2, \ldots$ can be obtained.
The idea is  to illustrate the behaviour of the product expansion by computing explicitly high order terms with  matrices in the $\mathrm{SU(2)}$
and the $\mathrm{SO(3)}$ Lie algebras.

\paragraph{Matrices $X$ and $Y$ in $\mathrm{SU(2)}$.}

In the first example we choose $A= i \,a \sigma_z$ and $B=i  \, \sigma_x$,
where $i = \sqrt{-1}$, $a$ is a real parameter and 
\begin{equation}\label{eq:sigma}
\sigma_x=\left(
\begin{array}{cc}
0&1\\1&0
\end{array} \right) , \quad
\sigma_y=\left(
\begin{array}{cc}
0&-i\\i&0
\end{array} \right), \quad
\sigma_z=\left(
\begin{array}{cc}
1&0\\0&-1
\end{array} \right) 
\end{equation}
are Pauli matrices. 
This instance is borrowed from  Quantum Mechanics, where $\exp[i(a\sigma_z+\varepsilon\sigma_x)]$ is a matrix that transforms the $\frac{1}{2}$-spin  wave function in a Hilbert space.

Using the scalar vector notation to write down a linear combination of Pauli matrices: 
$\vec v \cdot \vec \sigma = v_x\sigma_x+v_y\sigma_y+v_z\sigma_z$, the matrix exponential reads
\begin{equation}
\exp (i\vec v \cdot \vec \sigma)=\cos v\, I +i \,  \frac{\sin v}{v}\, \vec v \cdot \vec \sigma.
\end{equation}
In the sequel we work out the expansion 
\begin{equation}\label{eq:zx}
\e^{i(a\sigma_z+\varepsilon \sigma_x)}=\e^{ia\sigma_z}\,\e^{i \varepsilon W_1}\, \e^{i \varepsilon^2 W_2}\, \e^{i \varepsilon^3 W_3}\ldots 
\end{equation} 
up to order eleven in $\varepsilon$ and analyze the increasing accuracy of the product expansion as far as more terms are considered. The lhs in Eq. (\ref{eq:zx}) may be thought as a transformation involving $\sigma_x$ and $\sigma_z$. In turn, the rhs is a  pure $\sigma_z$ transformation, i.e.  $\exp(ia\sigma_z)$, followed by an infinite succession of transformations, $\exp(i \varepsilon^k W_k)$, whose effect should decrease with $k$.
The truncated product expansion is expected to be accurate as far as $\varepsilon  \ll a$.

\renewcommand*{\arraystretch}{2.}   
\begin{table}[t]
	\centering
	{\small	
		\begin{tabular}{ll}
			\hline \hline
			$k$ & $W_k(t)$ \\
			\hline
			1&$\dfrac{1}{2a}[S\sigma_x+(1-C)\sigma_y]$  \\
			2&$\dfrac{1}{4a^2}(2at-S)\sigma_z$  \\
			3& $\dfrac{1}{12a^3}\left\{ [6at+(C-4)S]\sigma_x-(1-C)^2\sigma_y \right\}$ \\
			4& $-\dfrac{1}{16a^4}[6at+(C-4)S]\sigma_z$ \\
			5& $\dfrac{1}{240a^5}\left\{ [56S-(4C+7)SC-10at(7+4C-2C^2)]\sigma_x +\right.$ \\
			\phantom{5} &$ \phantom{\dfrac{i}{240a^5}}  \left. 
			[4C^3-7C^2-28C+31+2at(2SC-4S+3at)]\sigma_y \right\}$ \\
			\hline \hline
		\end{tabular}
		\caption{First five orders in Bellman problem for generic $t$. The operators $C_k$ in Eq. (\ref{bell.0}) are obtained by taking $t=1$, i.e., $C_k=i  \, W_k(1)$. We have defined $S\equiv \sin(2at)$ and $C\equiv \cos(2at)$.} \label{tab:Bellman}
	}
\end{table}

In Table \ref{tab:Bellman} we write down explicitly  the first five contributions for a generic $t$ 
{(expressions for $k > 5$ are too involved to be collected here)}.
Wilcox--Bellman's formula (\ref{eq:zx}) corresponds then
to $t=1$. All the terms have been obtained with the recurrences of section \ref{section.3}, starting from
\begin{equation}\label{eq:dotW1}
\dot W_1(t)= \e^{-iat\sigma_z}\, \sigma_x \, \e^{iat\sigma_z}= C\sigma_x+S\sigma_y,
\end{equation}
where $C\equiv \cos(2at)$ and $S\equiv \sin(2at)$.

\begin{figure}[t]
	\centering
	\includegraphics[scale=0.45]{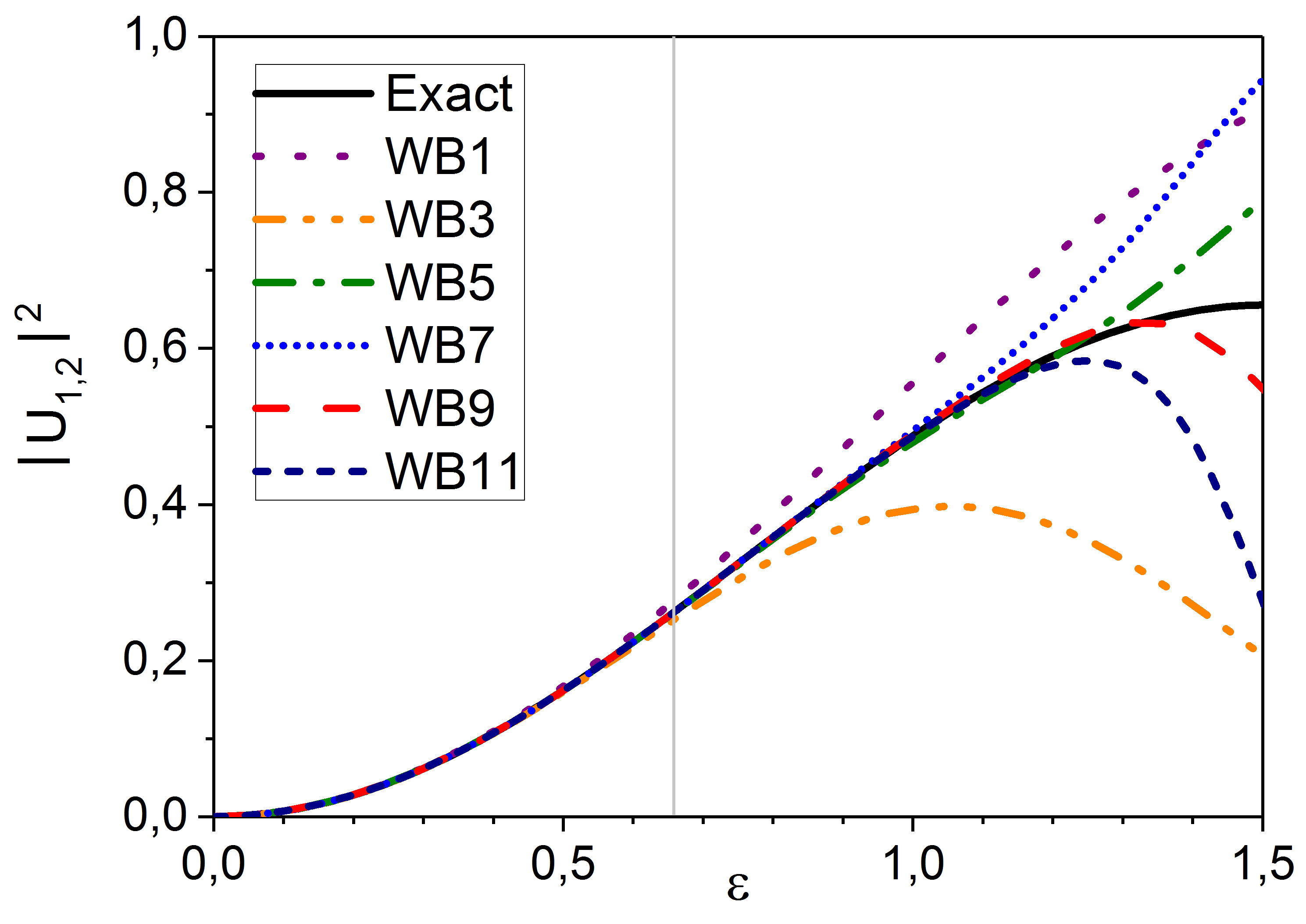}
	\caption[Wilcox--Bellman expansion example]{Accuracy of Wilcox--Bellman product expansion up to order eleven as a function of the 
		ratio $\varepsilon/a$, with $a=1$. The quantity plotted is the squared modulus of the non--diagonal element of the matrix. The vertical grey line stands for the convergence lower bound $\varepsilon=0.658$.}
	\label{fig:wb1-11}
\end{figure}

The formulas in Table \ref{tab:Bellman} show that $\varepsilon/a$ may be considered as an effective expansion parameter. In Figure \ref{fig:wb1-11} we  illustrate, for $a=1$,  the accuracy of the Wilcox--Bellman product expansion in the example at hand as a function of $\varepsilon$. We plot the squared modulus of the non--diagonal matrix element, say $|U_{1,2}|^2$, of Eq. (\ref{eq:zx}) for every {analytic} approximation up to order eleven
in $\varepsilon$, as well as the exact result. Even orders do not contribute in this test because $W_{2k}$ is always proportional to $\sigma_z$ and therefore $\exp(i\varepsilon^{2k}W_{2k})$ is a diagonal matrix.

As regards convergence of the product expansion, the lower bound of Eq. (\ref{conv.wil}) leads to
\begin{equation}
\int_0^1 \| \e^{-iat\sigma_z}\, \varepsilon\sigma_x \, \e^{iat\sigma_z}\|\,\mathrm{d}t= \varepsilon< 0.658.
\end{equation}
In turn, the behaviour of the curves in Figure \ref{fig:wb1-11} points out that convergence of the product expansion extends well
beyond that lower bound {for this particular example}.

\begin{figure}[t]
	\centering
	\includegraphics[scale=0.45]{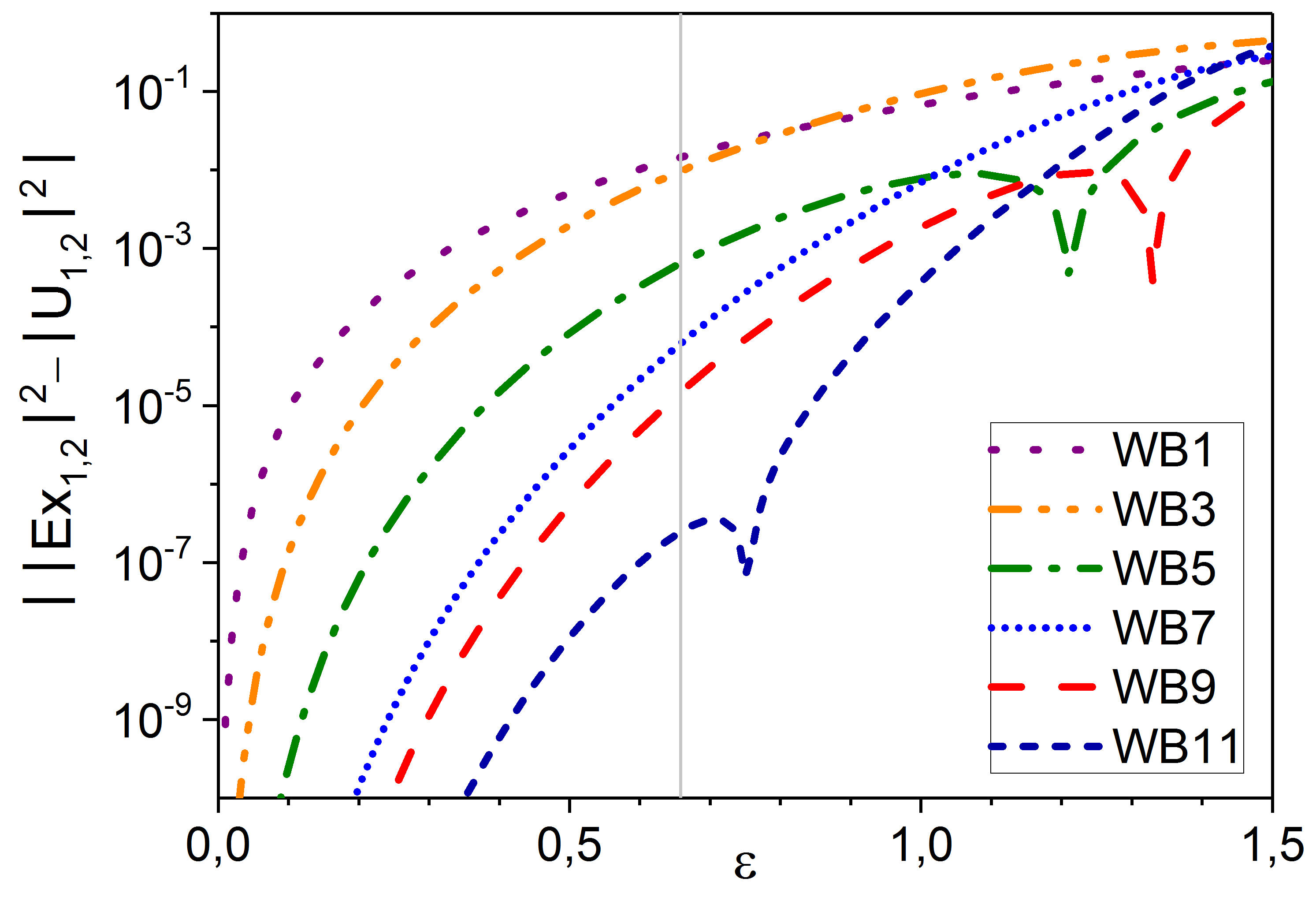}
	\caption{Absolute error of the approximations given by curves in Figure \ref{fig:wb1-11} with $a=1$. The vertical grey line is located at the value of the convergence lower bound $\varepsilon=0.658$.}
	\label{fig:WBerr1-11}
\end{figure}

Eventually, Figure \ref{fig:WBerr1-11} shows the logarithm of the absolute error in the approximations given by curves in Figure \ref{fig:wb1-11}. 


\paragraph{Matrices $X$ and $Y$ in $\mathrm{SO(3)}$.}

The second example refers to the matrix that describes a rotation in three dimensions defined by the vector $\vec{a}=\alpha \hat{a}$. Here $\alpha$ stands for the rotation angle around  the axis given by the unitary vector $\hat{a}$. 
A generic 3D rotation matrix can be written as $\exp(\vec{a}\cdot \vec{\rho})$, where the components of $\vec \rho$ are the three fundamental rotation matrices
\begin{equation}\label{eq:rot3D}
\rho_x=
\left(\begin{matrix}
0&0&0\\[-4mm]
0&0&-1\\[-4mm]
0&1&0
\end{matrix}\right), \qquad
\rho_y=
\left(\begin{matrix}
0&0&1\\[-4mm]
0&0&0\\[-4mm]
-1&0&0
\end{matrix}\right), \qquad
\rho_z=
\left(\begin{matrix}
0&-1&0\\[-4mm]
1&0&0\\[-4mm]
0&0&0
\end{matrix}\right).
\end{equation}
We study the particular case $\vec a \cdot \vec \rho= \alpha \big( \cos \theta \, \rho_z +\sin \theta \, \rho_x \big)$, {and compare the rotation of angle $\alpha$ around the unit vector $(\sin \theta,0,\cos\theta)$}
\begin{equation}
\mathrm{e}^{ \alpha \big( \cos \theta \, \rho_z +\sin \theta \, \rho_x \big)}
\end{equation}
with   {the sequence of transformations}
\begin{equation}\label{eq:zx2}
\mathrm{e}^{\alpha\cos\theta\,\rho_z}\,\mathrm{e}^{\alpha\sin\theta \,  W_1}\,\mathrm{e}^{(\alpha\sin\theta  )^2 W_2}\,\mathrm{e}^{(\alpha\sin\theta  )^3 W_3}\ldots  
\end{equation} 
In other words, the question we address is how {the} pure $z$--axis rotation {in Eq. (\ref{eq:zx2}), $\exp(\alpha\cos\theta\,\rho_z)$,} has to be corrected by an infinite composition of further rotations to reproduce the {one} defined by $\vec a \cdot \vec \rho$. When  $\alpha \sin\theta $ is small enough the approach is expected to converge since the expansion convergence lower bound  reads in this case $|\alpha \sin \theta| < 0.658$.

\begin{figure}[h!]
	\centering
	\includegraphics[scale=0.45]{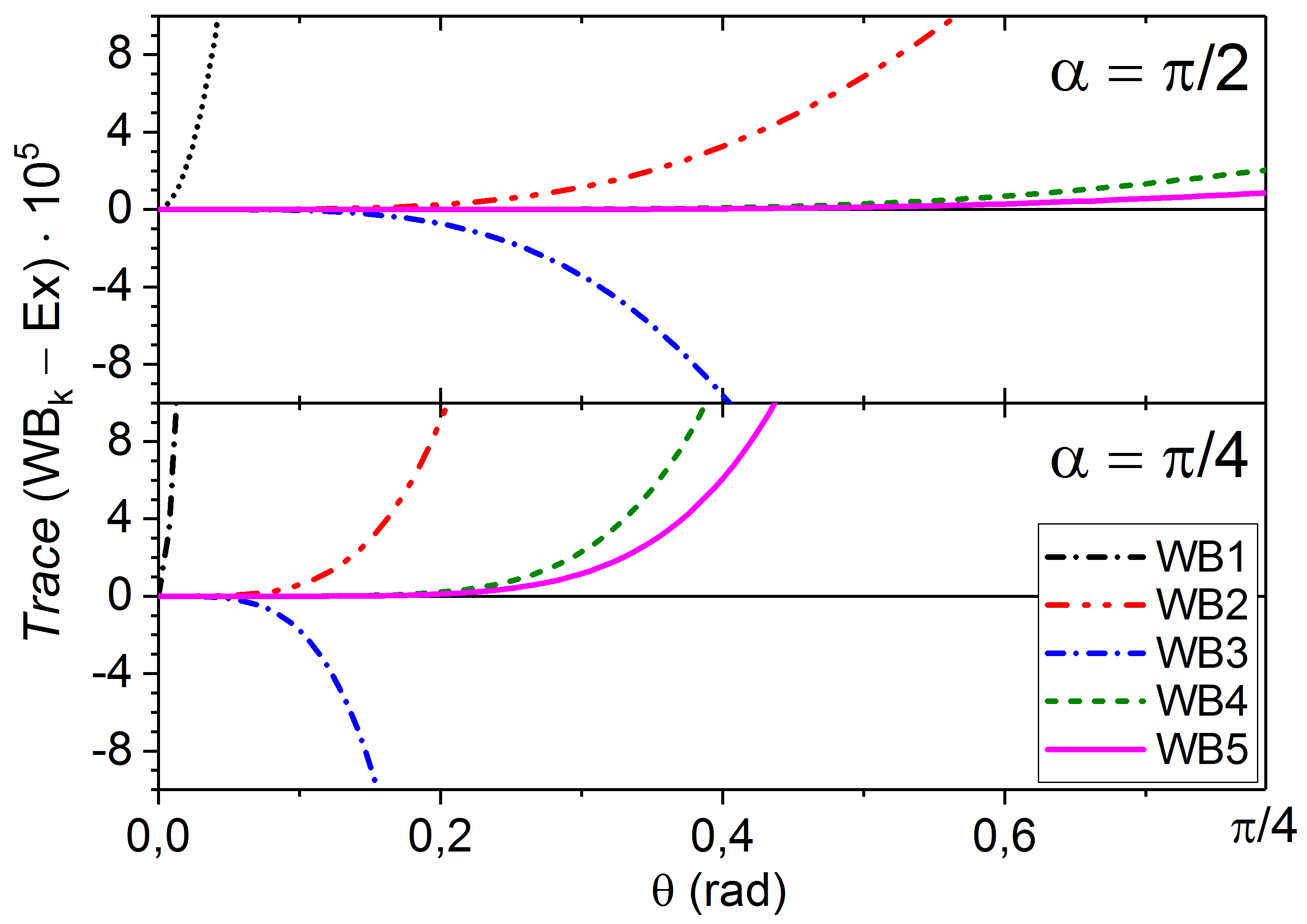}
	\caption{Error in the approximations to the matrix trace as a function of the rotation angle   $\theta$, for two values $\alpha=\pi/2$ and $\pi/4$.}
	\label{fig:WBrot}
\end{figure}

\begin{figure}[h!]
	\centering
	\includegraphics[scale=0.45]{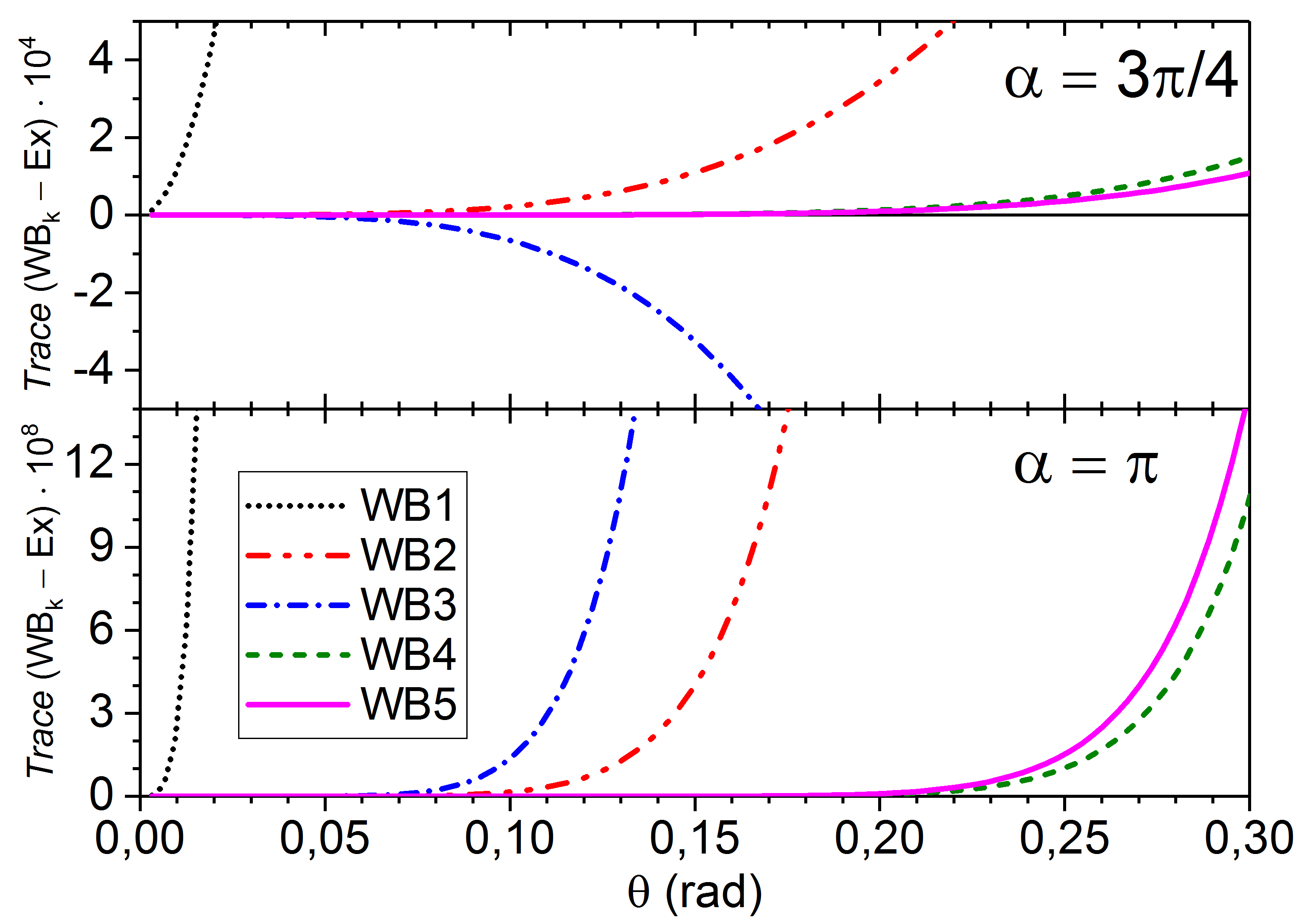}
	\caption{Error in the approximations to the matrix trace as a function of the rotation angle   $\theta$, for two values $\alpha=3\pi/4$ and $\pi$.}
	\label{fig:WBrot1}
\end{figure}

Here the accuracy of the product expansion will depend on both the rotation angle $\alpha$   and the relative orientation of the rotation axis, determined by the angle $\theta$. This is illustrated in Figures \ref{fig:WBrot} and \ref{fig:WBrot1} for the first five orders of approximation, with $\alpha=\pi/4,\pi/2,3\pi/4$ and $\pi$. 

In order to test the expansion, we have computed the matrix trace of the successive approximants and compared with the exact result: 
\begin{equation}
\mathrm{tr} \left(\exp \left(\alpha \big( \cos \theta \, \rho_z +\sin \theta \, \rho_x \big) \right) \right)=1 + 2\cos \alpha.   
\end{equation}    
The first two approximants  are simple enough to be written down:
\begin{equation}
\begin{aligned}
& \mathrm{tr}\left(\mathrm{e}^{\alpha\cos\theta\,\rho_z}\,\mathrm{e}^{\alpha\sin\theta \,  W_1} \right) = \Big(1+\cos (\alpha \cos \theta) \Big)\cos \left(2 \tan \theta \, \sin \frac{\alpha\cos\theta}{2}\right) \\
& \mathrm{tr}\left(\mathrm{e}^{\alpha\cos\theta\,\rho_z}\,\mathrm{e}^{\alpha\sin\theta \,  W_1}\,\mathrm{e}^{(\alpha\sin\theta  )^2 W_2}\right)= \cos\left(\frac{1}{2}\tan^2\theta \,(\sin(\alpha\cos\theta)-\alpha\cos \theta)\right)\cdot \\
& \phantom{ \mathrm{tr}\left(\mathrm{e}^{\alpha\cos\theta\,\rho_z}\,\mathrm{e}^{\alpha\sin\theta \,  W_1}\,\mathrm{e}^{(\alpha\sin\theta  )^2 W_2}\right) = }
\cos\left(\alpha \cos \theta -\frac{1}{2}\tan \theta \sin \frac{\alpha \cos \theta}{2}\right)
\end{aligned}
\end{equation}
Interestingly, the third approximation order is worse than the second one in all four cases. Also in the case $\alpha =\pi $ the fourth order is better than the fifth one.

\section{Conclusions}

When a linear system of differential equations, defined by the coefficient matrix $\lambda A(t)$, is transformed under $\exp(\lambda \int_0^t A(t') \mathrm{d}t' )$, the  coefficient matrix in the new representation becomes an infinite power series in $\lambda$, say 
$\tilde{A}(t)= \sum_{k \ge 1} \lambda^k \, a_k(t)$.
That is the first step of all matrix exponential methods to approximate the time-evolution operator. In the framework that we have introduced it is the first move in a sequence of exponential transformations that change the linear system from a representation to another with the goal that  dynamics will become less and less relevant. Choosing  the transformation $\exp(\lambda \int_0^t \tilde{A}(x) \mathrm{d}x )$ as the second move and iterating this procedure afterwards  yields the Fer expansion. Instead, choosing the transformation $\exp(\lambda \int_0^t a_1(x) \mathrm{d}x )$, i.e. the leading term of the new coefficient matrix,  opens up Wilcox expansion. The framework allows intermediate expansions  taking $\exp(\sum_{k=1}^n\lambda^k \int_0^t a_k(x) \mathrm{d}x)$ as initiator, as well as jumping between schemes, {according with the particular requirements of the problem at hand}.

{
	We have seen that the theory of linear transformations (or changes of picture in the language of Quantum Mechanics) provides a unified framework to deal
	with many different exponential perturbative expansions: whereas only one linear transformation reducing the dynamics to the trivial equation (\ref{trivial}) or
	to a system with a constant matrix
	renders the Magnus \cite{casas19cco} 
	and the Floquet-Magnus \cite{arnal20epe} expansion, respectively, a sequence of such transformations with different choices of the new matrices lead
	to Wilcox and Fer factorizations. From this perspective, other factorizations are possible depending on the particular problem at hand: one only has to 
	select appropriately the successive transformations.
}

In the case of Wilcox expansion we have provided an efficient recursive procedure to compute it. In addition, we have developed  a method to build up an explicit expression for any $W_n$ in terms of commutators. {This is possible by using similar tools as in the case of the Magnus expansion,
	namely by relating products of iterated integrals with the structure of the Hopf algebra of permutations, and by using special bases of nested
	commutators}.
A {sufficient} condition for the expansion convergence has also been obtained. 

We have presented some application examples of the results about Wilcox expansion. Firstly, we have shown how to obtain Zassenhaus formula from Wilcox expansion which, in turn, may be interpreted as its continuous analogue. Secondly, we point out that Wilcox expansion solves the problem of 
expanding  the exponential $\exp(A+\epsilon B)$ when $A$ and $B$ are non-commuting operators. We refer to it as Wilcox--Bellman expansion. Two practical cases on this respect have been analyzed up to high order. Interestingly, in one of them the convergence seems not to be uniform. 
For convenience, the interested reader can find in \cite{wilcoxcode} a \emph{Mathematica} code generating general
explicit expressions and recurrences for the Wilcox expansion.

{While a full assessment of the Wilcox expansion in comparison with Fer expansion is not the main purpose of this work, we can still mention here some of their
	most distinctive features. Both types of expansions construct the solution of Eq. (\ref{ccv.1}) as an infinite exponential factorization, but
	in Wilcox the exponent of each factor
	is proportional to successive powers of the expansion parameter $\lambda$, whereas in Fer each exponent contains an infinite sum of powers of $\lambda$.
	This means that, when truncated after a given number of transformations, say $n$, the Wilcox expansion differs from the exact solution in the power
	$\lambda^{n+1}$. In other words, each term $W_k(t)$
	in the Wilcox expansion collects the effect of the perturbation at order $k$. On the other hand, 
	the Fer expansion, when truncated after $n$ transformations, provides a much more accurate approximation. This is true, of course, if the infinite sums 
	involved in each transformation
	are exactly computed, an almost impossible task unless the time dependence of $A(t)$ is simple enough. 
	By contrast, we have explicit expressions
	for each exponent $W_k(t)$ in the Wilcox expansion for a generic $A(t)$ and, by using the same techniques as in the Magnus expansion, we can construct
	appropriate approximations to the iterated integrals if necessary. As the examples collected here and some other studies show \cite{zagury10ueo}, 
	Wilcox expansion can provide
	accurate results after only a few such transformations.
}

\vspace{1cm}

{\it The first three authors  has been funded by Ministerio de Ciencia e Innovaci\'on  (Spain) through projects MTM2016-77660-P and PID2019-104927GB-C21 (AEI/FE\-DER, UE),
	and by Universitat Jaume I (grants UJI-B2019-17 and GACUJI/2020/05). 
	The work of J.A.O.  has been partially supported by the Spanish MINECO (grant numbers AYA2016-81065-C2-2 and PID2019-109592GB-100).}


\end{document}